\documentclass[12pt]{article}
\usepackage{amsmath}
\newcommand{\be}{\begin{equation}}
\newcommand{\ee}{\end{equation}}
\newcommand{\ba}{\begin{eqnarray}}
\newcommand{\ea}{\end{eqnarray}}
\newcommand{\ban}{\begin{eqnarray*}}
\newcommand{\ean}{\end{eqnarray*}}

 \newcommand{\qed}{\hspace*{\fill}\rule{3mm}{3mm}\quad}
\newcommand{\Pf}{\noindent  {\em Proof.} }

\newcommand{\sect}[1]{\section{#1}  \setcounter{equation}{0}}

\begin{document}
\newtheorem{defn}[lem]{Definition}
\newtheorem{theo}[lem]{Theorem}
\newtheorem{cor}[lem]{Corollary}
\newtheorem{prop}[lem]{Proposition}
\newtheorem{rk}[lem]{Remark}
\newtheorem{ex}[lem]{Example}
\newtheorem{note}[lem]{Note}
\newtheorem{conj}[lem]{Conjecture}

\title{The Logarithmic Sobolev Inequality Along The Ricci Flow In Dimension 2}
\author{Rugang Ye \\ {\small Department  of Mathematics} \\
{\small University of California, Santa  Barbara}}
\date{July 20, 2007}
\maketitle

\noindent 1. Introduction \\
2. The main results \\
3. The Sobolev inequality and logarithmic Sobolev inequality on a Riemannian manifold of dimension 2 \\
4. Proofs of the main theorems \\

\sect{Introduction}

In [Y1], logarithmic Sobolev inequalities along the Ricci flow 
in dimensions $n \ge 3$ were obtained using Perelman's entropy monotonicity, which lead to Sobolev inequalities and $\kappa$-noncollpasing estimates.
At the time as [Y1] was posted, we had also obtained the corresponding results for the dimension $n=2$. (Indeed, the 
$C_{N,I}$-version of the 2-dimensional logarithmic Sobolev inequalities along the Ricci flow was obtained in 2004, the same time as 
the higher dimensional analogues in [Y1].) 
We did not post the paper on the 2-dimensional results  then, because we planned to expand the paper 
to include some further results. Now we have decided to post it without the planned expansion in order to clarify the issue. 

Note that the general theory in [Y1] covers all dimensions $n \ge 2$, except the logarithmic Sobolev 
inequality on a given Riemannian manifold which is used for dealing with the inital metric of the 
Ricci flow.  In [Y1], a logarithmic Sobolev inequality for a given Riemannian manifold of dimension $n \ge 3$ was 
presented. The $C_{N,I}$-version of the 2-dimensional case is similar, but the technical details are  different, which will be presented 
below.  On the other hand, the $\kappa$-noncollapsing estimate in dimension $n=2$ follows easily from 
the 3-dimensional result in [Y1]. This is seen by    
passing to the product with the circle.  We would like to point out that the logarithmic Sobolev inequality 
and the Sobolev inequality along the Ricci flow in dimension 2 also follow from the 3-dimensional results in 
[Y1] via the product with the unit circle. The approach presented in this paper is more direct. 

In the introduction section of [Y1] we explicitly stated that the results extended to the dimension 
$n=2$ and that the 2-dimensional case was going to be presented elsewhere.  It was entirely clear that 
the 2-dimensional case was not left open. One can also see this clearly from the presentation of the 
2-dimensional case in this paper.  Obviously, everything goes along the same lines as in [Y1]. 

\sect{The main results}

Consider a compact manifold $M$ of dimension $n=2$.  Let $g=g(t)$ be a smooth solution of the Ricci flow 
\ba
\frac{\partial g}{\partial t}=-2Ric
\ea
on $M \times [0, T)$ for some (finite or infinite) $T>0$ with a given initial metric 
$g(0)=g_0$.\\

\noindent {\bf Theorem A} {\it Let $A_1=A_1(g_0)$ and $A_2=A_2(g_0)$ be given by Theorem 
\ref{logtheorem2} or Theorem \ref{logtheorem4}. For each $\sigma>0$ and each $t \in [0, T)$  there holds
\ba \label{sobolevA}
\int_M u^2 \ln u^2 dvol &\le& \sigma \int_M (|\nabla u|^2 +\frac{R}{4} u^2)dvol  -\ln \sigma \nonumber 
\\ &&
+4(t+\sigma)A_1+A_2,
\ea
for all $u\in W^{1,2}(M)$ with  $\int_M u^2 dvol=1$, where
all geometric quantities  are associated with the metric $g(t)$ (e.g. the volume form $dvol$ and the scalar curvature $R$), except $A_1$ and $A_2$ which are associated with $g_0$.

Consequently, there holds for each $t \in [0, T)$
\ba \label{strongA}
\int_M u^2 \ln u^2 dvol \le  \ln \left[e^{1+A_1t+A_2}(\int_M (|\nabla u|^2 +\frac{R}{4} u^2)dvol+\frac{A_1}{4}) \right]
\ea
for all $u \in W^{1,2}(M)$ with $\int_M u^2 dvol=1$. 
} \\

The logarithmic Sobolev inequality in Theorem A is uniform for all time which lies below a given bound, but deteriorates as 
time becomes large. The next result takes care of large time under the assumption that the eigenvalue 
$\lambda_0$ of the initial metric is positive. \\

\noindent {\bf Theorem B} {\it Assume that the first eigenvalue $\lambda_0=\lambda_0(g_0)$ of the operator 
$-\Delta+\frac{R}{4}$ for  the initial metric $g_0$ is positive. Let $\delta_0=\delta_0(g_0)$ and $B_0=
B_0(g_0)$ be from Theorem \ref{logtheorem3} or Theorem \ref{logtheorem5}. Let $t \in [0, T)$ and $\sigma>0$ satisfy $t+\sigma \ge \frac{1}{4}\delta_0$.
Then there holds
\ba \label{sobolevB}
\int_M u^2 \ln u^2 dvol &\le& \sigma \int_M (|\nabla u|^2+\frac{R}{4}u^2)dvol -\ln \sigma +B_0
\ea
for all $u\in W^{1,2}(M)$ with  $\int_M u^2 dvol=1$, 
where all geometric quantities  are associated with the metric $g(t)$ (e.g. the volume form $dvol$ and the scalar curvature $R$), except the numbers 
$\sigma_0$  and $B_0$  which 
are associated with the initial metric $g_0$.

Consequently, there holds for each $t \in [0, T)$
\ba \label{strongB}
\int_M u^2 \ln u^2 dvol \le  \ln \left[e^{B_0+1}\int_M (|\nabla u|^2 +\frac{R}{4} u^2)dvol \right]
\ea
for all $u \in W^{1,2}(M)$ with $\int_M u^2 dvol=1$.
} \\

The above two results lead to Sobolev inequaltites along the Ricci flow. \\

\noindent {\bf Theorem C} \, \,  {Part I.} {\it  Assume $T<\infty$. For each $p>2$ there are positive constants $A$ and $B$ depending only on an upper bound for $\frac{1}{p-2}$, 
a nonpositive lower bound for $R_{g_0}$, a positive lower bound for $vol_{g_0}(M)$, an upper bound for $C_S(M,g_0)$,
and an upper bound for $T$, such that for each $t \in [0, T)$ and 
all $u \in W^{1,2}(M)$ there holds 
\ba \label{sobolevC1}
\left( \int_M |u|^{p} dvol \right)^{\frac{2}{p}} \le 
A\int_M (|\nabla u|^2+\frac{R}{4}u^2)dvol+B \int_M u^2 dvol,
\ea
where all geometric quantities except $A$ and $B$  are associated with $g(t)$. 
}\\ 
{Part II.} {\it Assume that $\lambda_0(g_0) > 0$. For each $p>2$  there is a  positive constant $A$ depending only on an upper bound for $\frac{1}{p-2}$, a nonpositive lower bound for $R_{g_0}$, a positive lower bound for $vol_{g_0}(M)$, an upper bound for $C_S(M,g_0)$,
and a positive lower bound for $\lambda_0(g_0)$, such that for each $t \in [0, T)$ and 
all $u \in W^{1,2}(M)$ there holds 
\ba \label{sobolevC2}
\left( \int_M |u|^{p} dvol \right)^{\frac{2}{p}} \le 
A\int_M (|\nabla u|^2+\frac{R}{4}u^2) dvol,
\ea
where all geometric quantities except $A$  are associated with $g(t)$. } 
\\

Next we state the  $\kappa$-noncollapsing estimate.  Let $ds^2$ denote the standard metric on the unit circle 
$S^1$. For the definition of the Sobolev constant $C_S$ see [Y1]. \\

\noindent {\bf Theorem D} Part I. \, \, {\it Assume that $T<\infty$.  There are 
positive constants $A$ and $B$ depending only on 
a nonpositive lower bound for $R_{g_0}$, a positive lower bound for $vol_{g_0}(M)$, an upper bound for $C_I(M,g_0)$ (or 
$C_{S}(M \times S^1, g_0 \times ds^2)$)
and an upper bound for $T$ with the following properties.
Let $L>0$ and $t \in [0, T)$. Consider the Riemannian manifold 
$(M, g)$ with $g=g(t)$. Assume $R\le \frac{1}{r^2}$ on a geodesic ball $B(x, r)$ with $0<r \le L$. Then 
there holds 
\ba \label{noncollapse1}
vol(B(x, r)) \ge \left(\frac{1}{A+BL^2}\right)^{\frac{3}{2}} r^2.
\ea
 } \\
Part II. {\it Assume that $\lambda_0(g_0) > 0$.  There is a positive constant $\alpha$ depending only on 
a nonpositive lower bound for $R_{g_0}$, a positive lower bound for $vol_{g_0}(M)$, an upper bound for $C_{N,I}(M,g_0)$
(or $C_{S}(M \times S^1, g_0 \times ds^2)$),
and a positive lower bound for $\lambda_0(g_0)$ with the following properies.
Let $t \in [0, T)$. Consider the Riemannian manifold 
$(M, g)$ with $g=g(t)$. Assume $R\le \frac{1}{r^2}$ on a geodesic ball $B(x, r)$ with $r>0$. Then 
there holds 
\ba \label{noncollapse2}
vol(B(x, r)) \ge \alpha r^2.
\ea
} \\

Next we address the question whether the above results can lead to uniform estimates independent of time without assuming 
$\lambda_0(g_0)>0$.   We have the following result in the case  that the genus of $M$ is zero. \\ 

\noindent {\bf Theorem E} {\it Assume that the genus of $M$ is zero, i.e. $M$ is diffeomorphic to 
the 2-sphere $S^2$ or the projective plane.   Then there holds  
$T \le  (4\pi)^{-1}vol_{g_0}(M)$.  Consequently,  Theorem A, Theorem B, Theorem C and Theorem D yield uniform estimates 
without any condition on time. (The estimates depend on an upper bound for $vol_{g_0}(M)$ in addition to 
the dependences stated in those theorems. To have simpler dependences we can rescale to achieve  $vol_{g_0}(M)=1$.) } \\

\noindent {\bf Remark}  This theorem contains two versions of results. One is with dependence  
on $C_{N,I}(M,g_0)$, i.e. when the constants $A_1, A_2, \delta_0$ and $B_0$ are from Theorem \ref{logtheorem2} 
and Theorem \ref{logtheorem3}.  The other is with dependence on $C_{S}(M \times S^1, g_0 \times ds^2)$,
i.e. when $A_1,A_2$ and $B_0$ are from Theorem \ref{logtheorem4} and Theorem \ref{logtheorem5}.
The uniform Sobolev inequality with dependence on $C_{N,I}(M,g_0)$ also follows from Hamilton's 
result on the monotonicity of the isoperimetric ratio [H1]. Our entropy approach provides a new perspective. 
Moreover, the uniform Sobolev inequality with dependence on $C_{S}(M \times S^1, g_0 \times ds^2)$ is stronger. \\

It seems that in the case of genus at least one our approach can only yield  estimates which 
depend on an upper bound on time. In particular, Theorem B is essentially not applicable for the following reason. 
If the genus is greater than one, then $\lambda_0(g_0)$ is always negative. If the  
genus equals  one, 
the eigenvalue $\lambda_0(g_0)$ is always negative when $g_0$ is nonflat, and is zero when 
$g$ is flat, see Theorem \ref{eigenvalue}.  
On the other hand, as an easy consequence of Hamilton's convergence theorems in [H2],  a uniform Sobolev inequality 
and a uniform logarithmic Sobolev inequality
hold true in the case of genus at least one. But the dependence on the intial metric is 
more involved. We state the Sobolev inequality.   
\\

\noindent {\bf Proposition F} {\it  Assume that the genus of $M$ is at least one. For each  
$p>2$ there are positive numbers 
$A$ and $B$ depending only on an upper bound for $\frac{1}{p-2}$, a positive lower bound for the volume $vol_{g_0}(M)$, an upper bound 
for the diameter $diam_{g_0}$, and an upper bound for the absolute value of the scalar curvature 
$|R_{g_0}|$ such that for each $t$ there holds 
\ba
\label{sobolevF}
\left( \int_M |u|^{p} dvol \right)^{\frac{2}{p}} \le 
A\int_M (|\nabla u|^2+\frac{R}{4}u^2)dvol+B \int_M u^2 dvol
\ea
for all $u \in W^{1,2}(M)$.  }

\sect{The Sobolev inequality and logarithmic Sobolev inequality on a Riemannian manifold of dimension 2}

Consider a closed Riemannian manifolds $(M, g)$ of dimension $n=2$.

\begin{prop} There holds for all $u \in W^{1,1}(M)$
\ba \label{poincare1}
\|u-\bar u\|_2 \le C_{N,I}(M, g) \|\nabla u\|_1
\ea
and 
\ba  \label{soboelv2}
\|u\|_2 \le C_{N,I}(M, g) \|\nabla u\|_1+ \frac{1}{vol_g(M)^{\frac{1}{2}}}\|u\|_1,
\ea
where $\bar u$ denotes the average of $u$ and $C_{N,I}(M,g)$ is the Neumann isoperimetric constant, see [Y1]. 
\end{prop}
\Pf It is well-known (see e.g. [L] or [Y2]) that 
\ba
\inf_{\sigma} \|u-\sigma\|_{\frac{n}{n-1}} \le C_{N,I}(M,g) \|\nabla u\|_1.
\ea
for a closed  $n$-dimensional Riemannian manifold $(M,g)$. Now we have $\frac{n}{n-1}=2$ because $n=2$. By minimizing the function 
$y(\sigma)=\int_M |u-\sigma|^2 dvol$ we see that its mimimum is achieved at 
$\sigma=\bar u$. This yields (\ref{poincare1}).  Then we have 
\ba
\|u\|_2 \le \|u-\bar u\|_2 +\|\bar u\|_2 \le C_{N,I}(M, g) \|\nabla u\|_1+ \frac{1}{vol_g(M)^{\frac{1}{2}}}\|u\|_1.
\ea
\qed \\

\begin{theo} \label{logtheorem1} There holds for all $u \in W^{1,1}(M)$ with $\int_M |u| dvol=1$
\ba \label{RLS1}
\int_M |u|\ln |u| dvol 
&\le& \ln \left(C^N_{I,g}(M) \|\nabla u\|_{1}
+\frac{1}{vol_g(M)^{\frac{1}{2}}}\right).
\ea
Consequently, there holds for all $u \in W^{1,2}(M)$ with $\int u^2 dvol =1$
\ba \label{RLS2}
\int_M u^2 \ln u^2 dvol &\le&  \ln \left( \int_M |\nabla u|^2 dvol 
+(C_{N,I}(M,g)+\frac{1}{vol_g(M)^{\frac{1}{2}}}) \right).
\ea
\end{theo}
\Pf  Assume $u \in W^{1,1}(M)$ with $\int_M |u| dvol =1$. 
Since $\ln$ is concave, we have by Jensen's inequality
\ba
\ln \int_M |u|^{2} dvol &=& \ln \int_M |u|\cdot  |u| dvol 
\ge \int_M |u| \ln |u| dvol 
\ea
It follows that 
\ba
\int_M |u|\ln |u| dvol &\le& \ln \int_M |u|^{2} dvol 
\nonumber \\
&\le&  \ln \left(C_{N,I}(M,g) \|\nabla u\|_{1}
+\frac{1}{vol_g(M)^{\frac{1}{2}}}\|u\|_1\right).
\ea

To show (\ref{RLS2}), we apply (\ref{RLS1}) to the function $u^2$ and then apply the 
inequality 
\ba
2|u \nabla u| \le \frac{1}{\epsilon} |\nabla u|^2 + {\epsilon} u^2
\ea
with the choice $\epsilon=C_{N,I}(M,g)$. 
\qed \\

 Applying Lemma 3.2 and Lemma 3.4 in [Y1] as in [Y1] 
we deduce from Theorem \ref{logtheorem1} the next two 
results.

\begin{theo}  \label{logtheorem2}  
 For each $\alpha>0$ and all $u \in W^{1,2}(M)$ with $\|u\|_2=1$ there holds for each $\sigma>0$
\ba \label{RLS3}
\int_M u^2 \ln u^2  dvol &\le& \sigma \int_M (|\nabla u|^2+\frac{R}{4}u^2)dvol -\ln \sigma
+A_1\sigma+A_2,
\ea
where 
\ba
A_1=A_1(g)=C_{N,I}(M, g)+vol_g(M)^{-\frac{1}{2}}-\frac{\min R^-}{4}
\ea
and $A_2=1$.
\end{theo}

\begin{theo} \label{logtheorem3} Assume that the first eigenvalue $\lambda_0=\lambda_0(g)$ of the operator 
$-\Delta+\frac{R}{4}$ is positive.  For each $\sigma \ge \delta_0$ and all 
$u \in W^{1,2}(M)$ with $\|u\|_2=1$ there holds 
\ba \label{RLS4}
\int_M u^2 \ln u^2 \le \sigma \int_M(|\nabla u|^2+\frac{R}{4}u^2)
-\ln \sigma+ B_0,
\ea
where
\ba \label{delta-0}
\delta_0=\delta_0(g)=
\left(\lambda_0(g)+C_{N,I}(M,g)+vol_g(M)^{-\frac{1}{2}}
\frac{\min R^-}{4}\right)^{-1}
\ea 
and
\ba \label{B-0}
B_0=B_0(g)=\ln(1+\lambda_0(g)^{-1}(C_{N,I}(M,g)+vol_g(M)^{-\frac{1}{2}}-\frac{\min R^-}{4}))-1.
\ea
\end{theo}

An alternative way of deriving a suitable logarithmic Sobolev inequality 
for a Riemannian manifold of dimension 2 is in terms of the  Sobolev constant $C_{S}(M \times S^1, g_0 \times ds^2)$.
Since this constant can be estimated in terms of $C_{N,I}(M, g_0)$, the following two theorems are 
stronger than the above two theorems.  To keep this paper more streamlined, their proofs are 
presented in [Y2]. 

\begin{theo} \label{logtheorem4} There are positive constants $A_1=A_1(g)$ and $A_2=A_2(g)$ depending only on 
a positive lower bound for $vol_{g_0}(M)$, a nonpositive lower bound for 
$R_{g_0}$, and an upper bound for $C_{S}(M \times S^1, g_0 \times ds^2)$ such that 
\ba
\label{RLS5}
\int_M u^2 \ln u^2 dvol  \le \sigma \int_M |\nabla u|^2 dvol -\ln \sigma+ A_1 \sigma+A_2
\ea
for each $\sigma>0$ and all $u \in W^{1,2}(M)$ with $\int_M u^2 dvol=1$. 
\end{theo}

\begin{theo} \label{logtheorem5} Assume that the first eigenvalue $\lambda_0=\lambda_0(g)$ of the operator 
$-\Delta+\frac{R}{4}$ is positive. There are constants $\delta_0>0$ and 
$B_0$ depending only on a positive lower bound for $vol_{g_0}(M)$, a nonpositive lower bound for 
$R_{g_0}$, an upper bound for $C_{S}(M \times S^1, g_0 \times ds^2)$, and 
a positive lower bound for $\lambda_0(g_0)$ such that  
\ba \label{RLS6}
\int_M u^2 \ln u^2 \le \sigma \int_M(|\nabla u|^2+\frac{R}{4}u^2)
-\ln \sigma+ B_0 
\ea
for each $\sigma \ge \delta_0$ and all $u \in W^{1,2}(M)$ with $\int_M u^2 dvol=1$. 
\end{theo}

Finally, we present a result on the eigenvalue $\lambda_0$ in the case that the genus of 
$M$ equals $1$, i.e. $M$ is diffeomorphic to the 2-torus or the Klein bottle. 

\begin{theo} \label{eigenvalue} Assume that the genus of $M$ equals one.  Then $\lambda_0(g)<0$ for 
all $g$ on $M$ except when $g$ is a flat metric, in which case $\lambda_0(g)=0$. 
\end{theo}
\Pf Let $g$ be a metric on $M$.  For $u\equiv 1$ we have 
\ba
\int_M (|\nabla u|^2+\frac{R}{4})dvol=\int_M \frac{R}{4} dvol=0.
\ea
Hence $\lambda_0(g) \le 0$.  On the other hand, if $R \equiv 0$, then  we obviously have $\lambda_0(g) \ge 0$ 
and hence  $\lambda_0(g)=0$. 

Next let a metric $g_0$ satisfy $\lambda_0(g_0)=0$. Let $g=g(t)$ be the smooth solution 
of the Ricci flow with $g(0)=g_0$ on its maximal time interval $[0, T)$. (By [H], we have $T=\infty$, but we don't need this fact.)  Fix $t_1 \in (0, T)$. Let $u_1$ be a positive eigenfunction for the eigenvalue $\lambda_0(g(t_1))$ 
associated with the metric $g(t_1)$, such that 
$\int_M u_1^2 dvol=1$ with respect to $g(t_1)$.  Let $f=f(t)$ be the smooth solution of the 
equation 
\ba \label{f}
\frac{\partial f}{\partial t}=-\Delta f+|\nabla f|^2-R
\ea
on $[0, t_1]$ with $f(t_1)=-2\ln u_1$.  Note that (\ref{f}) is equivalent to
\ba
\frac{\partial v}{\partial t}=-\Delta v+Rv,
\ea
where $v=e^{-f}$. So the solution $f(t)$ exists. We also infer 
$\frac{d}{dt}\int_M v dvol=0$, and hence $\int_M v dvol =1$ for all $t \in [0, t_1]$. 
 
 We set $u=e^{-\frac{f}{2}}.$  By [P, (1.4)] we then have
\ba
\frac{d}{dt} \int_M (|\nabla u|^2+\frac{R}{4}u^2)dvol 
=\frac{1}{4} \frac{d}{dt}\int_M (|\nabla f|^2+R)e^{-f}dvol \ge \frac{1}{2} \int_M |Ric+\nabla^2 f|^2 e^{-f}dvol.
\nonumber \\
\ea
It follows that 
\ba
0 &\ge& \lambda_0(g(t_1)) 
\ge \lambda_0(g_0) + \frac{1}{2} \int_0^{t_1} \int_M |Ric+\nabla^2 f|^2 e^{-f}dvol dt
\nonumber \\
&=& \frac{1}{2} \int_0^{t_1} \int_M |Ric+\nabla^2 f|^2 e^{-f}dvol dt.
\ea
We deduce 
\ba
Ric+\nabla^2 f=0
\ea
on $[0, t_1]$. Hence $g_0$ is a steady Ricci soliton. By [CK, Proposition 5.20], it is Ricci flat.  \qed \\

\sect{Proofs of the main theorems}

\noindent {\bf Proof of Theorem A and Theorem B} These two theorems follow from 
Theorem 4.2 in [Y1] together with Theorem \ref{logtheorem2}, Theorem \ref{logtheorem3},
Theorem \ref{logtheorem4} and Theorem \ref{logtheorem5}. \qed \\

\noindent {\bf Proof of Theorem C} Part I.  For a given $p>2$, we set $\mu=
\frac{2p}{p-2}$. Then $\mu>2$ and $p=\frac{2\mu}{\mu-2}$. For $0<\sigma \le 1$ we derive from
(\ref{sobolevA})
that for each $t \in [0, T)$  there holds
\ba \label{sobolevAproof}
\int_M u^2 \ln u^2 dvol &\le& \sigma \int_M (|\nabla u|^2 +\frac{R}{4} u^2)dvol  -\frac{\mu}{2}\ln \sigma \nonumber 
\\ &&
+4(T+1)A_1+A_2,
\ea
for all $u\in W^{1,2}(M)$ with  $\int_M u^2 dvol=1$. Applying Theorem 5.3 in [Y1] with 
$\sigma^*=1$ and $\Psi=\frac{R}{4}$ we then infer 
\ba 
\label{Dsobolev3}
\|u\|^2_{p} \le  \left(\frac{1}{4}\right)^{1-\frac{2}{\mu}} c
\left(\int_M (|\nabla u|^2+\frac{R}{4})dvol+(4-\frac{\min R^-}{4}) \int_M u^2 dvol\right)
\ea
for all $u \in W^{1,2}(M)$, where $c=c(\bar C, \frac{1}{\mu-2})=c(\bar C, 
\frac{4}{p-2})$ and 
\ba \label{barC1}
\bar C=2^{\frac{2}{p-2}} e^{\frac{p}{2(p-2)}-\frac{3}{16}\min \frac{R^-}{4}+2(T+1)A_1+\frac{1}{2}A_2}. 
\ea
 Since $\min R^- \ge \min R^-_{g_0}$, we arrive at the 
desired Sobolev inequality (\ref{sobolevC1}). \\
Part II. This is similar to Part I. For $p>2$, $\mu=\frac{2p}{p-2}$ and $0<\sigma \le 1$ we 
derive from
(\ref{sobolevB})
that for each $t \in [0, T)$  there holds
\ba \label{sobolevBproof}
\int_M u^2 \ln u^2 dvol &\le& \sigma \int_M (|\nabla u|^2 +\frac{R}{4} u^2)dvol  -\frac{\mu}{2}\ln \sigma 
+B_0
\ea
for all $u\in W^{1,2}(M)$ with  $\int_M u^2 dvol=1$. As in Part I we then arrive at (\ref{Dsobolev3}) with
\ba \label{barC2}
\bar C=2^{\frac{2}{p-2}} e^{\frac{p}{2(p-2)}-\frac{3}{16}\min \frac{R^-}{4}+\frac{1}{2}B_0}. 
\ea
Since $\lambda_0$ is nondecreasing along the Ricci flow, we have 
\ba \label{lambda}
\int_M u^2 dvol \le \frac{1}{\lambda_0(g_0)} \int_M (|\nabla u|^2+\frac{R}{4}u^2) dvol.
\ea
Combining  (\ref{Dsobolev3}) with (\ref{lambda}) and the inequality 
$\min R^- \ge \min R^-_{g_0}$ we arrive at (\ref{sobolevC2}). \qed \\

\noindent {\bf Proof of Theorem D} Part I.  Obviously, the product metric $g^*(t) =
g(t) \times ds^2$ is a smooth solution of the Ricci flow on 
$(M \times S^1) \times [0, T)$ with the initial metric $g^*(0)=g_0 \times ds^2$. 
Let $L>0$ and $t \in [0, T)$. Consider the Riemannian manifold 
$(M, g)$ with $g=g(t)$. Assume $R \le \frac{1}{r^2}$ on a geodesic ball 
$B(x, r)$ with $0<r\le L$. Then $R_{g^*} \le \frac{1}{r^2}$ on 
$B(x, r) \times B_{S^1}(s_0, r)$, where $g^*=g^*(t)$, $s_0$ is an arbitary point in 
$S^1$, and $B_{S^1}(s_0, r)$ denotes the geodesic ball of center $s_0$ and 
radius $r$ in $S^1$.  Let $B_{M \times S^1}((x, s_0), r)$ denotes 
the corresponding geodesic ball in $(M\times S^1, g^*)$. It is easy to see that 
$B_{M \times S^1}((x, s_0), r) \subset B(x,r) \times B_{S^1}(s_0, r)$. Hence 
we infer $R_{g^*} \le \frac{1}{r^2}$ on $B_{M \times S^1}((x, s_0), r)$. 
By Theorem $\mbox{E}^*$ in [Y1] we then have 
\ba
vol(B(x,r)) \cdot vol(B(s_0, r)) \ge vol(B_{M \times S^1}((x, s_0), r))\ge \left(\frac{1}{2^6 A+2BL^2}\right)^{\frac{3}{2}} r^3,
\ea
where $A$ and $B$ depend only on 
a nonpositive lower bound for $R_{g_0}$, a positive lower bound for $vol_{g_0}(M)$, an upper bound for  
$C_{S}(M \times S^1, g_0 \times ds^2)$,
and an upper bound for $T$.
But $vol(B(s_0, r))\le 2\min\{r, \pi\}$. Hence we arrive at 
\ba\label{volume1}
vol(B(x,r))  \ge  \left(\frac{1}{2^6 A+2BL^2}\right)^{\frac{3}{2}} \frac{r}{2\min\{r, \pi\}} r^2,
\ea 
which implies (\ref{noncollapse1}) with redefined $A$ and $B$.  \\
Part II This is similar to Part I. We apply Theorem E instead of 
Theorem $\mbox{E}^*$ in [Y1] and arrive at 
\ba  \label{volume2}
vol(B(x,r))\ge \left(\frac{1}{2^{6}A}\right)^{\frac{3}{2}} \frac{r}{2\min\{r, \pi\}} r^2,
\ea
where $A$ depends only on a nonpositive lower bound for $R_{g_0}$, a positive lower bound for $vol_{g_0}(M)$, an upper bound for $C_{N,I}(M,g_0)$
(or $C_{S}(M \times S^1, g_0 \times ds^2)$),
and a positive lower bound for $\lambda_0(g_0 \times ds^2)$. But $\lambda_0(g_0 \times ds^2) 
\ge \lambda_0(g_0)$. Indeed, we have for $u \in W^{1,2}(M \times S^1)$ 
\ba
&& \int_{M \times S^1} (|\nabla_{g^*(0)}  u|_{g^*(0)}^2+ \frac{R_{g^*(0)}}{4}u^2) dvol_{g^*(0)}  \ge  
\int_{S^1} ds \int_M  (|\nabla_{g_0} u|_{g_0}^2 + \frac{R_{g_0}}{4}u^2) dvol_{g_0} 
\nonumber \\ && \ge 
\int_{S^1} \lambda_0(g_0) \int_M u^2 dvol_{g_0} ds = \lambda_0(g_0) \int_{M \times S^1} 
u^2 dvol_{g^*(0)}.
\ea  
Hence we arrive at the desired estimate (\ref{noncollapse2}), where 
\ba
\alpha=\frac{1}{2} \left(\frac{1}{2^{6}A}\right)^{\frac{3}{2}}.
\ea
\qed \\

\noindent {\bf Proof of Theorem E} We have 
\ba
\frac{d}{dt} vol_{g(t)}(M) = -\int_M R dvol=-4\pi.
\ea
Hence 
\ba
vol_{g(t)}(M)=vol_{g_0}(M)-4\pi t
\ea
and then 
\ba
t \le (4\pi)^{-1} vol_{g(0)}(M).
\ea
\qed

\end{document}